\documentclass{amsart}
\usepackage{amssymb} 
\newtheorem{thm}{Theorem}[section]
\newtheorem{cor}[thm]{Corollary}
\newtheorem{lem}[thm]{Lemma}
\newtheorem{prop}[thm]{Proposition}
\theoremstyle{definition}

\theoremstyle{remark}
\newtheorem{rem}{Remark}

\numberwithin{equation}{section}

\begin{document}

\title[Left 3-Engel elements in groups]{Left 3-Engel elements in groups}%
\author{Alireza Abdollahi}%
\address{Department of Mathematics, University of Isfahan, Isfahan 81746-73441, Iran}%
\email{a.abdollahi@sci.ui.ac.ir}%

\thanks{This research  supported by Isfahan University Grant no. 801031.
}%
\begin{abstract}
In this paper we study left 3-Engel elements in groups. In
particular, we prove that for any prime $p$ and any left 3-Engel
element $x$ of finite $p$-power order in a group $G$, $x^p$ is in
the Baer radical of $G$. Also it is proved that $\left<x,y\right>$
is nilpotent of class $4$ for every two left 3-Engel elements in a
group $G$.
\end{abstract}
\maketitle
\section{Introduction and results}
Let $G$ be any group and $n$ be a non-negative integer. For any
two elements $a$ and $b$ of $G$ we define inductively  $[a, _n
b]$ the $n$-Engel commutator of the pair $(a,b)$,  as follows:
$$[a,_0b]:=a, \;\; [a,b]:=a^{-1}b^{-1}ab, \;\;\text{and}\;\;
[a,_n b]:=[[a, _{n-1} b],b] \;\;\text{for all}\;\; n>0.$$
 An element $x$ of $G$ is called a left $n$-Engel
element if $[g,_n x]=1$ for all $g\in G$. We denote by $L_n(G)$,
the set of all left $n$-Engel elements of $G$. So $L_0(G)=1$,
$L_1(G)=Z(G)$ the centre of $G$, and it can be easily seen that
$$L_2(G)=\{x\in G \;|\; \left<x\right>^G \;\text{is abelian}\;
\},$$ where $\left<x\right>^G$ denotes the normal closure of $x$
in $G$. Therefore $L_2(G)$ is contained in $B(G)$ the Baer radical
of $G$, and in particular it is contained in $HP(G)$ the
Hirsch-Plotkin radical of $G$. In general for an arbitrary  group
$K$ it is not necessary that  $L_n(K)\subseteq HP(K)$. For
suppose, for a contradiction, that $L_n(K)\subseteq HP(K)$ for all
$n$ and all groups $K$. By a deep result of Ivanov \cite{Iv},
there is a finitely generated  infinite group $M$ of exponent
$2^k$ for some  positive integer $k$. Suppose that $k$ is the
least integer with this property, so every finitely generated
group of exponent dividing $2^{k-1}$ is  finite. Let $a$ be any
element of $M$ of order $2$ and $x$ an arbitrary element of $M$.
It is easy to see that $[x,_m a]=[x,a]^{(-2)^{m-1}}$ for all
positive integers $m$. Thus every element of order $2$ of $M$ is
in $L_{k+1}(M)$. So by hypothesis, $M/HP(M)$ is of exponent
dividing $2^{k-1}$ and so it is finite. Since $M$ is finitely
generated, $HP(G)$ is also.  But this yields that $HP(G)$ is a
periodic finitely generated nilpotent group and so  it is finite.
It follows that $M$ is finite, a contradiction.\\
 Hence  the
question which naturally arises, is that: What is the least
positive integer $n$ such that $L_n(G) \nsubseteq HP(G)$? To
study this question it should first study the case $n=3$, since
$$L_0(G)\subseteq L_1(G)\subseteq L_2(G)\subseteq\cdots\subseteq
L_n(G)\subseteq\cdots .$$ The main object of this paper is to
study $L_3(G)$. As far as we know there is no example of a group
$G$ for which $L_3(G)\nsubseteq HP(G)$. The corresponding subset
 to $L_n(G)$  which can be similarly defined, is
$R_n(G)$ the set of all right $n$-Engel elements of $G$: an
element $x$ of $G$ is called a right $n$-Engel element if $[x,_n
g]=1$ for all $g\in G$. There is a well-known relation between
these two subsets of $G$ due to Heineken \cite[Theorem
7.11]{Rob2}: $(R_n(G))^{-1}\subseteq L_{n+1}(G)$ for all integers
$n>0$. It is clear that  $R_0(G)=1$, $R_1(G)=Z(G)$ and by a
result of Kappe \cite[Corollary 1 Theorem 7.13]{Rob2}, $R_2(G)$ is
a characteristic subgroup of $G$. It is known also that
$R_2(G)\subseteq L_2(G)$ \cite[Theorem 7.13 (i)]{Rob2}. Recently
Newell \cite{Ne} has shown that the normal closure of every
element of $R_3(G)$ is nilpotent of class at most $3$. This
shows, of course, that $R_3(G)\subseteq B(G)\leq HP(G)$. An early
example was given by Macdonald \cite{Mac} shows that the inverse
and square of a right 3-Engel element need not be right 3-Engel,
thus $R_3(G)$ is not necessarily a subgroup. But as we  show in
Corollary \ref{cor1}, $L_3(G)$ is closed under the taking all
powers. It is interesting to note that $(R_3(G))^2\subseteq
L_3(G)$ and $(R_3(G))^4\subseteq R_3(G)$ \cite[Remark to Theorem
1]{Ne}.\\

The main results of this paper are the followings.
\begin{thm}\label{thm1}
Let $G$ be a group and $p$ be a prime number. If $x\in L_3(G)$
and  $x^{p^n}=1$ for some integer $n>1$, then
$\left<x^p\right>^G$ is soluble of derived length at most $n-1$
and $x^p$ belongs to $B(G)$, the Baer radical of $G$. In
particular,  $x^p$ belongs to $HP(G)$, the Hirsch-Plotkin radical
of $G$.
\end{thm}

\begin{thm} \label{thm2}
Let $G$ be any group and $a,b\in L_3(G)$. Then $\left<a,b\right>$
is nilpotent of class at most $4$.
\end{thm}
\section{proofs}

\begin{lem}\label{lem1}
Suppose that $G$ is an arbitrary group and $x,y\in G$. Then
$[y,_3 x]=[y^{-1},_3x]=1$ if and only if $\left<x,x^y\right>$ is
nilpotent of class at most $2$.
\end{lem}
\begin{proof}
Since
$$[y^{\epsilon},_3x]=[x^{-1},[x^{-1},[x^{-1},y^{\epsilon}]]]^{x^3}$$
where $\epsilon\in\{-1,1\}$, we have
\begin{align*}
  &[y,_3x]=[y^{-1},_3x]=1   \\
& \Leftrightarrow
 [x^{-1},[x^{-1},[x^{-1},y]]]=[x^{-1},[x^{-1},[x^{-1},y^{-1}]]]=1\\
&\Leftrightarrow
[x^{-x^{-y}},x^{-1}]=[x^{-x^{-y^{-1}}},x^{-1}]=1\\
&\Leftrightarrow [x^{-y},_2x^{-1}]=[x^{-y^{-1}},_2 x^{-1}]=1\\
&\Leftrightarrow [x^{-y},_2 x^{-1}]=[x^{-1},_2 x^{-y}]=1\\
&\Leftrightarrow \left<x^{-1},x^{-y}\right>=\left<x,x^y\right> \;
\text{is nilpotent of class at most}\; 2.
\end{align*}
\end{proof}
So we have the following characterization of left 3-Engel
elements in a group. We denote by $\mathcal{N}_2$ the class of
nilpotent groups of class at most $2$.
\begin{cor}\label{cor1}
For an arbitrary group $G$,
$$L_3(G)=\{x\in G \;|\; \left<x,x^y\right>\in\mathcal{N}_2
\;\text{for all}\; y\in G\}.$$ In particular every power of a
left 3-Engel element is also a left $3$-Engel element.
\end{cor}
\begin{prop}
A group generated by a set of left $3$-Engel elements of finite
orders such that their  orders are pairwise coprime is abelian.
\end{prop}
\begin{proof}
It is enough to show that $[a,b]=1$ for any two left 3-Engel
elements $a$ and $b$ such that $\gcd(|a|,|b|)=1$. By Corollary
\ref{cor1}, $K=\left<a,a^b\right>$ and $H=\left<b,b^a\right>$ are
 both nilpotent. Thus $[a,b]=a^{-1}a^b=(b^{-1})^ab\in K\cap H$. Since
$\gcd(|a|,|b|)=1$ and $H$ and $K$ are  both nilpotent,
$\gcd(|H|,|K|)=1$. It follows that $[a,b]=1$, as required.
\end{proof}
\begin{lem}\label{lem11}
Let $p$ be a prime number, $G$ be a group and $x\in L_3(G)$. If
$x^{p^n}=1$ for some integer $n\geq 2$, then $x^{p^{n-1}}\in
L_2(G)$.
\end{lem}
\begin{proof}
Let $y$ be an arbitrary element of $G$. By Corollary \ref{cor1},
$\left<x,x^y\right>$ is nilpotent of class at most $2$. Thus
\begin{align*}
&[(x^{-y})^{p^{n-1}},x^{p^{n-1}}]=\\
&[(x^{-y})^{p^{n-2}},x^{p^{n-1}}]^p=  \;\;\; \text{since}\;\; [x^{-y},x] \in Z(\left<x,x^y\right>)\\
&[(x^{-y})^{p^{n-2}},x^{p^n}]=1
\end{align*}
But
$[y,x^{p^{n-1}},x^{p^{n-1}}]=[(x^{-y})^{p^{n-1}},x^{p^{n-1}}]=1$.
This completes the proof.
\end{proof}

{\bf Proof of Theorem \ref{thm1}.}
 By Lemma \ref{lem11} and Corollary
\ref{cor1}
$$1\trianglelefteq \left<x^{p^{n-1}}\right>^G\trianglelefteq
 \left<x^{p^{n-2}}\right>^G\trianglelefteq
\cdots\trianglelefteq \left<x^{p}\right>^G$$ is a series of
normal subgroups of $G$ with abelian factors. This implies that
$K=\left<x^p\right>^G$ is soluble of derived length at most
$n-1$. By Corollary \ref{cor1}, $x^p$ and so all its conjugates
in $G$  belong to  $L_3(G)$ and in particular they are in
$L_3(K)$.
 Now a result of Gruenberg
 \cite[Theorem 7.35]{Rob2} implies that $B(K)=K$. But
 $K\trianglelefteq G$ which yields  that $x^p\in K\leq B(G)\leq
 HP(G)$. \;\;\;$\Box$\\

In the following calculations, one must be careful with notation.
As usual $u^{g_1+g_2}$ is shorthand notation for $u^{g_1}u^{g_2}$.
This means that
$$u^{(g_1+g_2)(h_1+h_2)}=u^{(g_1+g_2)h_1}u^{(g_1+g_2)h_2}$$ which
does not have to be equal to $u^{g_1(h_1+h_2)}u^{g_2(h_1+h_2)}$.
We also have that $$u^{(g_1+g_2)(-h)}=((u^{g_1}u^{g_2})^{-1})^h$$
which is equal to $u^{-g_2h-g_1h}$. This does not have to be the
same as $u^{-g_1h-g_2h}$.\\
 The following remark easily follows from Corollary
\ref{cor1}.
 We use in sequel this remark, sometimes without any reference.
\begin{rem}
Let $G$ be any group and $a\in L_3(G)$. Then
$$[a,x]^{a^2}=[a,x]^{2a-1} \;\text{and}\;\;
[a,x]^{a^{-1}}=[a,x]^{-a+2} \;\text{for all}\; x\in G.$$
\end{rem}

\begin{lem}\label{lem2}
Let $G$ be any group and $a,b\in L_3(G)$. Then
$$[\left<a,b\right>,\left<a,b\right>]=\left<[a,b],[a,b]^a,[a,b]^b,[a,b]^{ab}\right>.$$
\end{lem}
\begin{proof}
It is enough to show   that
$N=\left<[a,b],[a,b]^a,[a,b]^b,[a,b]^{ab}\right>$ is a normal
subgroup of $\left<a,b\right>$. We have
\begin{align*}
&[a,b]^{a^2}=[a,b]^{2a-1}\in N\\
&[a,b]^{ba}=[a,b]^{ab[b,a]}=[a,b]^{1+ab-1}\in N\\
&[a,b]^{aba}=[a,b]^{a^2b[b,a]}=[a,b]^{1+a^2b-1}
&=[a,b]^{1+(2a-1)b-1}=[a,b]^{1+2ab-b-1}\in N
\end{align*}
Thus $N^a\leq N$ (I). Also we can write
\begin{align*}
&[a,b]^{a^{-1}}=[a,b]^{-a+2}\in N \\
&[a,b]^{ba^{-1}}=[a,b]^{a^{-1}b[b,a^{-1}]}=[a,b]^{-a+a^{-1}b-a^{-1}}=[a,b]^{-a-ab+2b+a^{-1}}\in
N\\
&[a,b]^{aba^{-1}}=[a,b]^{b[b,a^{-1}]}=[a,b]^{-a^{-1}+b-a^{-1}}\in
N
\end{align*}
 Therefore $N^{a^{-1}}\leq N$ (II). It follows from (I) and (II)
 that $N=N^a$. Now since $[a,b]^{ab}=[a,b]^{-1+ba+1}\in N$ we have
 $M=\left<[b,a],[b,a]^b, [b,a]^a,[b,a]^{ba}\right>$ is equal to
 $N$. Thus by the symmetry between $a$ and $b$, we have $M^b=M$
 and so $N^b=N$. Hence $N\trianglelefteq \left<a,b\right>$. This
 completes the proof.
\end{proof}
\begin{lem}\label{lem3}
Let $G$ be any group and $a,b\in L_3(G)$. Then both of the
commutators $[a,b]$ and $[a,b]^{ab}$ belong to
$C_G([a,b]^a,[a,b]^b)$.
\end{lem}
\begin{proof}
By Lemma \ref{lem1}, $\left<a,a^b\right>\in \mathcal{N}_2$ and
$\left<b,b^a\right>\in\mathcal{N}_2$, so $[[a,b],[a,b]^a]=1$ (I)
and $ [[a,b],[a,b]^b]=1$ (II). It follows from (I) that
$[[a,b]^b,[a,b]^{ab}]=1$ and (II) yields that
$[[a,b]^a,[a,b]^{ba}]=1$ (III). But $[a,b]^{ba}=[a,b]^{1+ab-1}$
so (I) and (III) imply that $[[a,b]^a,[a,b]^{ab}]=1$. This
completes the proof.
\end{proof}

{\bf Proof of Theorem \ref{thm2}.}
 We first prove that the derived subgroup of $\left<a,b\right>$ is
 abelian. Since $\left<a,b\right>'$ is the normal closure of
 $[a,b]$ in $\left<a,b\right>$, it is enough to show that $[a,b]$
 is in the centre of $\left<a,b\right>'$ and by Lemma \ref{lem3},
it suffices to prove   that $[[a,b],[a,b]^{ab}]=1$.\\
Let $a^b=a_1$, $a^{b^2}=a_2$ and let $c=[a,b]$. Then
$\left<a,a_1\right>$, $\left<a_1,a_2\right>$ and
$\left<a_1,a_2\right>$ are each nilpotent of class at most 2, by
Lemma \ref{lem1}.\\
Now $c=a^{-1}a_1$ and $c^b=a_1^{-1}a_2$ and these commute by Lemma
\ref{lem1}. But
\begin{align*}
1=[c^b,c]&=[a_1^{-1}a_2,a^{-1}a_1]\\
&=[a_1^{-1}a_2,a_1][a_1^{-1}a_2,a^{-1}]^{a_1}\\
1&=[a_2,a_1][a_1,a]^{a_2a_1}[a_2,a]^{-a_1}\\
1&=[a_2,a_1][a_1,a][a,a_2]
\end{align*}

Hence $[a_1,a,a_2]=1$. Now $[a_1,a]=[a,b,a]=[c,a]$ and
$a_2=ac^2[c,b]$. Since $[c,a]$ commutes with $a$ and $c$, it
follows that $[[c,a],[c,b]]=1$.  Now It follows from  Lemma
\ref{lem3} that $[c^a,c^b]=1$. Now by Remark 1 we have
$c^{b^{-1}}=c^{-b}c^2$. It follows that
$[c^{ab},c]=[c^a,c^{b^{-1}}]=1$.\\
 Thus $\left<a,b\right>$ is metabelian. Now for completing the
proof, it is enough to show that $[a,b,x_1,x_2,x_3]=1$ (*) for all
$x_1,x_2,x_3\in\{a,b\}$. Since $\left<a,b\right>$ is metabelian we
have
$[a,b,x_1,x_2,x_3]=[a,b,x_{\sigma(1)},x_{\sigma(2)},x_{\sigma(3)}]$
for all permutations $\sigma$ on the set $\{1,2,3\}$.  Now $a,b\in
L_3(G)$, it is easy to see that the equality (*) is satisfied   by
all $x_1,x_2,x_3\in\{a,b\}$. This completes the
proof.  \;\;\;$\Box$\\

The author could not prove a similar result to Theorem \ref{thm2}
for groups generated by $3$ left 3-Engel elements.
 We finish
this paper by the following question.\\

\noindent {\bf Question.} Is there a function
$f:\mathbb{N}\rightarrow \mathbb{N}$ such that every nilpotent
group generated by $d$ left
$3$-Engel elements is nilpotent of class at most $f(d)$?\\

\noindent {\bf Acknowledgement.} The author is  indebted to the
referee who has read this paper  carefully and patiently and
pointed out a serious error in the
origin proof of Theorem \ref{thm2} and \\


\end{document}